\documentclass[12pt,amsfonts]{article}
\usepackage{graphicx}
\usepackage{latexsym}
\usepackage{amssymb}
\usepackage{amsmath}
\usepackage{enumerate}

\usepackage{layout}
\usepackage{eufrak}
\newtheorem{prop}{Proposition}
\newtheorem{lemma}{Lemma}

\newtheorem{corollary}{Corollary}
\newtheorem{theorem}{Theorem}
\newtheorem{remark}{Remark}
\newtheorem{example}{Example}

\def\real{{\mathord{{\rm I\kern-2.8pt R}}}}        
\def\inte{{\mathord{{\rm I\kern-2.8pt N}}}}

\def\sZZ{{\rm Z\kern-2.8ptem{}Z}}

\def\z{{\mathchoice
  {\sZZ}
  {\sZZ}
  {\rm Z\kern-0.30em{}Z}
  {\rm Z\kern-0.25em{}Z} }}
\def\sQQ{{\kern 0.27em \vrule height1.45ex width0.03em depth0em
          \kern-0.30em \rm Q}}
\def\qu{{\mathchoice
    {\sQQ}
    {\sQQ}
  {\kern 0.225em \vrule height1.05ex width0.025em depth0em \kern-0.25em \rm Q}
  {\kern 0.180em \vrule height0.78ex width0.020em depth0em \kern-0.20em \rm Q}
        }}
\def\sCC{{\kern 0.27em \vrule height1.45ex width0.03em depth0em
          \kern-0.30em \rm C}}
\def\complex{{\mathchoice
    {\sCC}
    {\sCC}
  {\kern 0.225em \vrule height1.05ex width0.025em depth0em \kern-0.25em \rm C}
  {\kern 0.180em \vrule height0.78ex width0.020em depth0em \kern-0.20em \rm C}
        }}

\newcommand{\ba}{\begin{array}}
\newcommand{\ea}{\end{array}}
\newcommand{\be}{\begin{equation}}
\newcommand{\ee}{\end{equation}}
\newcommand{\bea}{\begin{eqnarray}}
\newcommand{\eea}{\end{eqnarray}}
\newcommand{\beaa}{\begin{eqnarray*}}
\newcommand{\eeaa}{\end{eqnarray*}}

\def\z{\zeta}

\font\tenmath=msbm10 \font\sevenmath=msbm7 \font\fivemath=msbm5
\newfam\mathfam \textfont\mathfam=\tenmath
\scriptfont\mathfam=\sevenmath \scriptscriptfont\mathfam=\fivemath

\def \={{\buildrel {\rm (law)} \over =}}

\def\qed{ \hfill \vrule width.25cm height.25cm depth0cm\smallskip}

\newcommand{\basa}{\begin{assumption}}
\newcommand{\easa}{\end{assumption}}

\newcommand{\bas}{\begin{assum}}
\newcommand{\eas}{\end{assum}}

\newcommand{\ignore}[1]{}
\begin{document}

\renewcommand{\thefootnote}{\fnsymbol{footnote}}

\renewcommand{\thefootnote}{\fnsymbol{footnote}}

\title{The determinant of the iterated  Malliavin matrix and the density of a couple of multiple integrals }

\author{   {David Nualart} \thanks{Supported by the NSF grant DMS-1208625.}    \\
  Department of Mathematics \\
  University of Kansas\\
   Lawrence, Kansas 66045, USA \\
   nualart@ku.edu
\and
\\
{Ciprian  A. Tudor} \thanks{  Supported by the CNCS grant PN-II-ID-PCCE-2011-2-0015. }\\ Laboratoire Paul Painlev\'e \\ Universit\'e de Lille I \\ F-59655 Villeneuve d'Ascq, France \\and 
Academy of Economical Studies, Bucharest, Romania\\
tudor@math.univ-lille1.fr
}
\date{}

 \maketitle
 \abstract
 {The aim of this paper is to show an estimate for the determinant of the covariance of a two-dimensional vector of multiple stochastic integrals of the same order in terms of a linear combination of the expectation of the determinant of its iterated Malliavin matrices. As an application we show that  the vector is absolutely continuous  if and only if its components are proportional.}

\vskip0.3cm

{\bf 2010 AMS Classification Numbers:}   60F05, 60H05, 91G70.

 \vskip0.3cm


\section{Introduction}
A basic result in Malliavin calculus says that if the Malliavin matrix $\Lambda=(\langle DF_i ,DF_j \rangle_H)_{1\le i,j \le d}$ of a $d$-dimensional random vector  $\mathbf{F}=(F_1, \dots, F_d)$ is nonsingular almost surely, then this vector has an absolutely continuous law with respect to the Lebesgue measure in $\mathbb{R}^d$. In the special case of vectors whose components belong to a finite sum of  Wiener chaos, Nourdin, Nualart and Poly proved in \cite{NoNuPo} that the following conditions are equivalent:
\begin{itemize}
\item[(a)] The law of  $\mathbf{F}$ is not absolutely continuous.
\item[(b)] $E\det \Lambda =0$.
\end{itemize}
A natural question  is the relation between  $E\det \Lambda$ and the determinant of the covariance matrix  $C$ of the random vector $\mathbf{F}$. Clearly if $\det C=0$, then the  components of $\mathbf{F}$ are linearly dependent and the law of $\mathbf{F}$ is not absolutely  continuous, which implies $E\det \Lambda=0$. The converse is not true if $d\ge 3$. For instance,  the vector $(F_1, F_2, F_1F_2)$, where $F_1$ and $F_2$ are two non-zero independent random variables  in the first chaos, satisfies  $\det \Lambda  =0$ but $\det C\not =0$.

The purpose of  this paper is to show the equivalence between  $E\det \Lambda=0$ and $\det C=0$ in the particular case of a  two-dimensional random vector  $(F,G)$ whose components are multiple stochastic integrals of the same order $n$.  This implies that the random vector $(F,G)$  has an absolutely continuous law with respect to the Lebesgue measure on $\mathbb{R}^2$ if and only if its components are proportional, as in the Gaussian case.
This  result was established for $n=2$ in   \cite{NoNuPo}, and for $n=3,4$ in \cite{T}. Our proof in the general case $n\ge 2$ is based on the   notion of iterated Malliavin matrix and the computation of the expectation of its determinant.

In connection with this equivalence we will derive  an inequality relating  $E\det \Lambda$ and $\det C$, which has its own interest. In the case of double stochastic integrals, that is if $n =2$,  it was proved  in \cite{NoNuPo}  that
\[
E \det \Lambda \ge 4\det C.
\]
We extend this inequality proving that 
\[
E \det \Lambda \ge c_n\det C
\]
holds for $n=3,4$ with $c_3= \frac 94$ and $c_4= \frac{16} 9$.  For $n\ge 5$ we obtain a more involved inequality, where in the left hand side we have a linear combination (with positive coefficients) of the expectation of the iterated Mallavin matrices of $(F,G)$ of order $k$ for
$1\le k\le \left [\frac {n-1} 2\right]$ (see Theorem \ref{thm1} below).  

The paper is organized as follows. In Section 2 we present some preliminary results and notation. Section 3 contains a  general decomposition of the determinant of the iterated Malliavin matrix of a two-dimensional random vector into a sum of squares.  In Section 4 we prove our main result which is based on a further decomposition of the determinant of the iterated Malliavin matrix of a vector whose components are multiple stochastic integrals. Finally, the application to the characterization of absolutely continuity is obtained in Section 5.

\section{Preliminaries}

We briefly describe the tools from the analysis on Wiener space that we will need in our work. For complete presentations, we refer to \cite{N} or \cite{NPbook}.
Let $H$ be a real and separable Hilbert space and consider   an isonormal process $(W(h), h\in H)$. That is, $(W(h), h\in H)$ is a Gaussian family of centered   random variables on a probability space $(\Omega, {\cal{F}}, P)$  such that $E W(h) W(g)= \langle f, g\rangle _{H}$ for every $h, g\in H$. Assume that the $\sigma$-algebra ${\cal{F}}$ is generated by $W$. 

For  any integer $n\geq 1$ we denote  by ${\cal{H}} _{n}$ the $n$th  Wiener chaos generated by $W$. That is, ${\cal{H}} _{n}$ is the vector subspace of $L^{2}(\Omega)$ generated by  the random variables $\left( H_{n} (W(h)), h\in H, \Vert h\Vert _H=1\right)$ where $H_{n}$ the Hermite polynomial of degree $n$.   We denote by $H_0$ the space of constant random variables. Let $H^{\otimes n}$ and $H^{\odot n}$ denote, respectively, the $n$th tensor product and the $n$th symmetric tensor product of $H$.
For any $n\geq 1$, the mapping $I_{n}(h^{\otimes n})= H_{n}(W(h))$ can be extended to an isometry between the symmetric tensor product space  $H^{\odot n}$ endowed with the norm $\sqrt{n!} \Vert \cdot \Vert _{H ^{\otimes n}} $ and the $n$th  Wiener chaos ${\cal{H}}_{n}$.  For any $f\in H^{\odot n}$, the random variable $I_{n}(f)$ is called the multiple Wiener It\^o integral of $f$ with respect to $W$.

Consider $(e_{j})_{j\geq 1}$ a complete orthonormal system in $H$ and let $f\in H^{\odot n}$, $ g\in H^{\odot m}$ be two symmetric tensors with $n,m\geq 1$.  Then
\begin{equation}\label{f}
f= \sum_{j_{1}, \dots, j_{n} \geq 1} f _{j_{1},\dots, j_{n}} e_{j_{1}} \otimes \cdots \otimes e_{j_{n}}
\end{equation}
and
\begin{equation}\label{g}
g= \sum_{k_{1}, \dots , k_{m}\geq 1} g _{ k_{1},\dots , k_{m}} e_{k_{1}}\otimes \cdots \otimes e_{k_{m}},
\end{equation}
where the coefficients are given by  $f _{j_{1}, \dots , j_{n}}  = \langle f, e_{j_{1}} \otimes \cdots \otimes e_{j_{n}}\rangle$ and $g _{ k_{1},\dots , k_{m}} =\langle g, e_{k_{1}}\otimes \cdots\otimes e_{k_{m}} \rangle $. 
  These coefficients are symmetric, that is,  they satisfy
  $f_{ j_{\sigma (1)} ,\dots, j_{\sigma (n)}}= f _{j_{1},\dots , j_{n}}$ and $ g _{ k_{\pi(1)},\dots , k_{\pi(m)}}=g_{ k_{1},\dots, k_{m}}$ for every permutation $\sigma $ of the set $\{1,\dots, n\}$ and for every permutation $\pi$ of the set $\{1,\dots, m\}$.   
  
  Note that, throughout the paper we will usually omit   the subindex  $H^{\otimes k}$ in the notation for the norm and the 
  scalar product in $H^{\otimes k}$ for any $k\ge 1$. 

If $f\in H^{\odot n}$, $ g\in H^{\odot m}$ are symmetric   tensors given by (\ref{f}) and (\ref{g}), respectively, then the contraction of order  $r$  of $ f$ and $g$ is given by
\begin{eqnarray}
f\otimes _{r}g &=& \sum_{ i_{1},\dots, i_{r} \geq 1} \sum_{j_{1}, \dots ,j_{n-r}\geq 1} \sum_{k_{1},\dots, k_{m-r}\geq 1} f _{i_{1},\dots, i_{r} , j_{1}, \dots , j_{n-r}}g _{i_{1},\dots , i_{r} , k_{1},  \dots, k_{m-r}}\nonumber\\
&&\times \left( e_{j_{1}}\otimes  \cdots \otimes e_{j_{n-r}}\right) \otimes \left( e_{k_{1}}\otimes \cdots\otimes e_{k_{m-r}}\right) \label{contra}
\end{eqnarray}
for every $r=0, \dots , m\wedge n$. In particular $f\otimes _{0}g= f\otimes g$.
Note that $f\otimes _{r} g $ belongs to $H^{\otimes (m+n-2r)}$ for every $r=0,\dots, m\wedge n
$ and it is not in general symmetric. We will denote by $f\tilde{\otimes }_{r}g$ the symmetrization of $f\otimes _{r}g$.  In the particular case when $H= L^{2}(T, {\cal{B}}, \mu)$ where $\mu$ is a sigma-finite measure without atoms, (\ref{contra}) becomes
\begin{eqnarray}
&&(f\otimes _{r}g) (t_{1},\dots, t_{m+n-2r})\nonumber 
 =  \int_{T^{r}}d\mu(u_{1}) \cdots d\mu (u_{r}) f(u_{1},  \dots, u_{r}, t_{1},\dots, t_{n-r}) \\
 && \quad \times g(u_{1},\dots,u_{r}, t_{n-r+1},\dots , t_{m+n-2r}).  \label{contra2}
\end{eqnarray}

An important role will be played by  the following product formula for multiple Wiener-It\^o integrals: if $f\in H^{\odot n}$, $ g\in H^{\odot m}$ are symmetric tensors, then
\begin{equation}
\label{prod}
I_{n}(f) I_{m}(g)= \sum_{r=0} ^{m\wedge n} r! C_{m}^{r} C_{n}^{r} I_{m+n-2r} \left(f\tilde{\otimes }_{r}g\right).
\end{equation}

We will need some elements of the Malliavin calculus with respect to the isonormal Gaussian process $W$.  Let $\mathcal{S}$ be the set of all smooth and cylindrical random variables of
the form
\begin{equation}
F=\varphi \left( W(h_{1}),\dots , W(h _{n})\right) ,  \label{v3}
\end{equation}%
where $n\geq 1$, $\varphi:\mathbb{R}^{n}\rightarrow \mathbb{R}$ is a infinitely
differentiable function with compact support, and $h_{i}\in  H$ for $i=1,..,n$.
If $F$ is given by (\ref{v3}), the {\it Malliavin derivative} of $F$ with respect to $W$ is the element of $%
L^{2}(\Omega ; H)$ defined as
\begin{equation*}
DF = \sum_{i=1}^{n}\frac{\partial \varphi}{\partial x_{i}}\left( W(h
_{1}),\ldots ,W(h _{n})\right)h _{i}.
\end{equation*}
By iteration, one can
define the $k$th derivative $D^{(k)}F$ for every $k\geq 2$, which is an element of $L^{2}(\Omega ;%
  H^{\odot  k})$.
For $k\geq 1$, ${\mathbb{D}}^{k,2}$ denotes the closure of $%
\mathcal{S}$ with respect to the norm $\Vert \cdot \Vert_{\mathbb{D}^{k,2}}$, defined by
the relation
\begin{equation*}
\Vert F\Vert _{\mathbb{D}^{k,2}}^{2}\;=\;E\left[ |F|^{2}\right] +\sum_{i=1}^{k}E\left(
\Vert D^{(i)}F\Vert _{ H^{\otimes i}}^{2}\right) .
\end{equation*}
If $ F=I_n(f)$, where $f\in H^{ \odot n}$  and $ I_{n}(f)$ denotes the multiple integral of order $n$ with respect to $W$, then 
\begin{equation*}
DI_{n}(f) =n \sum_{ j=1} ^\infty   I_{n-1} \left(  f \otimes _1 e_{j} \right) e_j.
\end{equation*}
More generally, for any $1 \le k\le n$, the iterated Malliavin derivative of $I_n(f)$ is given by
\begin{equation*}
D^{(k)}I_{n}(f) = \frac {n!} {(n-k)!}  \sum_{j_1, \dots, j_k \ge 1}     I_{n-k} \left(  f_{j_1, \dots, j_k} \right) e_{j_{1}} \otimes  \cdots \otimes e_{j_{k}},
\end{equation*}
where
\begin{equation}\label{fjj}
f_{j_1, \dots, j_k} = f \otimes _k  (e_{j_{1}} \otimes  \cdots \otimes e_{j_{k}}).
\end{equation}
We denote by $\delta $ the adjoint of the operator $D$, also called the
{\it divergence operator} or {\it Skorohod integral}. A random element $u\in L^{2}(\Omega ;H)$ belongs to the domain of $\delta $, denoted $\mathrm{Dom}\delta $, if and
only if it verifies
\begin{equation*}
\left |E \langle DF,u\rangle _{H} \right|\leq c_{u}\,\sqrt{E(F^2)}
\end{equation*}%
for any $F\in \mathbb{D}^{1,2}$, where $c_{u}$ is a constant depending only
on $u$. If $u\in \mathrm{Dom}\delta $, then the random variable $\delta (u)$
is defined by the duality relationship 
\[
E(F\delta (u))=E\langle DF,u\rangle _{ H},  
\]
which holds for every $F\in {\mathbb{D}}^{1,2}$. 
If $F=I_n(f)$ is a multiple stochastic integral of order $n$, with $f\in H^{\odot n}$, then $DF$ belongs to the domain of $\delta$ and 
\begin{equation} \label{delta}
\delta DF =nF.
\end{equation}

\section{Decomposition of the  determinant of the iterated Malliavin matrix}
In this section we obtain a decomposition into a sum of a squares for the determinant of the iterated Malliavin matrix  of a $2$-dimensional random vector.  We recall that if   $F,G$ are two random variables  in the space $\mathbb{D}^{1,2}$,    \textit{the Malliavin matrix} of the random vector $(F,G)$ is the defined as the following $2\times2$  random matrix
\[
\Lambda =\left(  \begin{array}{cc}
\Vert DF\Vert _{H} ^{2}&\langle DF, DG\rangle_{H}\\ \langle DF, DG\rangle_{H} &\Vert DF\Vert _{H} ^{2}
\end{array} \right) .
\]
More generally, fix $k\ge 2$ and suppose that $F,G$ are two random variables in $\mathbb{D}^{k,2}$. The   \textit{$k$th iterated Malliavin matrix of   the vector $(F,G)$} is defined as  
\[
\Lambda ^{(k)} =\left(
   \begin{matrix}   
      \Vert D^{(k)}F \Vert  ^{2}_{H^{\otimes k}} & \langle D^{(k)}F, D^{(k)} G \rangle _{ H ^{\otimes k} } \\
    \langle D^{(k)}F, D^{(k)} G \rangle _{ H ^{\otimes k} } &    \Vert D^{(k)}G \Vert  ^{2}_{H^{\otimes k}}  \\
   \end{matrix}
   \right),
   \]
  We set $\Lambda ^{(1)} =\Lambda$. 
For every $j_{1},\dots ,j_{k}\geq 1$, we will write
\[
D^{(k)}_{j_1, \dots, j_k} F=  \langle   D^{(k)} F , e_{j_{1}} \otimes  \cdots \otimes e_{j_{k}} \rangle _{ H ^{\otimes k} }.
\]
 The next proposition provides an expression of the determinant of   the iterated Malliavin matrix of a random vector as a sum of squared random variables.
\begin{prop}
Suppose that $(F, G)$ is a $
2$-dimensional random vector whose components belong to $\mathbb{D}^{k,2}$ for some $k\ge 1$. Let $\Lambda ^{(k)}$ be the $k$th iterated Malliavin matrix of $(F,G)$.  Then
\begin{equation} \label{e1}
\det \Lambda ^{(k)} =  \frac{1}{2}\sum_{ i_{1},\dots, i_{k}, l_{1},  \dots, l_{k}\geq 1}\left(  D^{(k)} _{i_1, \dots, i_k}  F  D^{(k)} _{l_1, \dots,l_k} G -  D^{(k)}  _{l_1, \dots, l_k} F  D^{(k)}  _{i_1, \dots, i_k}  G \right) ^{2}.
\end{equation}
\end{prop}
{\bf Proof: } For every $k \geq 1$ we  have
\[
\Vert D^{(k)}F \Vert _{H^{\otimes k}} ^{2} = \sum _{i_{1} , \dots, i_{k}\geq 1}   \left(D^{(k)}_{i_1, \dots, i_k}  F\right)^2 ,
\]
\[
\Vert D^{(k)}G \Vert _{H^{\otimes k}} ^{2} = \sum _{i_{1} , \dots, i_{k}\geq 1}   \left(  D^{(k)}_{i_1, \dots, i_k} G \right)^2
\]
and
\[
 \langle D^{(k)}F, D^{(k)} G \rangle _{ H ^{\otimes k} }=\sum _{i_{1},\dots, i_{k}\geq 1}  D^{(k)}_{i_1, \dots, i_k}  FD^{(k)}_{i_1, \dots, i_k} G.
 \] 
 Thus
\begin{eqnarray*}
\det \Lambda ^{(k)} &=& \sum _{i_{1} , \dots, i_{k}\geq 1}   \left(D^{(k)}_{i_1, \dots, i_k}  F\right)^2
 \sum _{i_{1} , \dots, i_{k}\geq 1}   \left(  D^{(k)}_{i_1, \dots, i_k} G \right)^2 \\
 && -\left( \sum _{i_{1},\dots, i_{k}\geq 1}  D^{(k)}_{i_1, \dots, i_k}  FD^{(k)}_{i_1, \dots, i_k} G \right) ^{2}  \\
&=&\frac{1}{2}\sum_{ i_{1},.., i_{k}, l_{1}, .., l_{k}\geq 1}\left(  D^{(k)} _{i_1, \dots, i_k}  F  D^{(k)} _{l_1, \dots,l_k} G -  D^{(k)}  _{l_1, \dots, l_k} F  D^{(k)}  _{i_1, \dots, i_k}  G \right) ^{2}.
\end{eqnarray*}
\qed

\section{The iterated Malliavin matrix of a two-dimen\-sional vector of multiple integrals}

Throughout this section, we assume  that  the components of the random vector $(F,G)$ are multiple Wiener-It\^o integrals. More precisely,   we will fix $n,m\geq 1$ and we will  consider the vector
$$(F,G)= (I_{n}(f), I_{m}(g) )$$
where $f\in H^{\odot n}$ and $g\in H^{\odot m}$.
Since for every $1\leq k\leq   \min(n,m)$, 
\[
D_{i_1, \dots, i_k}  ^{(k)} F  =\frac{n!}{(n-k) !} I_{n-k}\left( f _{i_{1},\dots, i_{k}}\right)
\]
(with $f _{i_{1},\dots, i_{k}}$ defined by (\ref{fjj})) and 
\[
D^{(k)} _{i_1, \dots, i_k}  G  =\frac{m!}{(m-k) !} I_{m-k}\left( g _{i_{1},\dots, i_{k}}\right)
\]
formula (\ref{e1}) reduces to 
\begin{eqnarray*}
\det \Lambda ^{(k)} &=&\frac{1}{2} \left(\frac{n!}{(n-k)!} \frac{m!}{(m-k)! } \right) ^{2} \sum_{ i_{1},\dots, i_{k}, l_{1}, \dots, l_{k}\geq 1}
\left[ I_{n-k} \left( f _{i_{1},\dots, i_{k}}\right)I_{m-k} \left( g _{l_{1},\dots, l_{k}}\right) \right. \\
&& \left. -I_{n-k} \left( f _{l_{1},\dots, l_{k}}\right)I_{m-k} \left( g _{i_{1},\dots, i_{k}}\right)\right] ^{2}.
\end{eqnarray*}
By the product formula for multiple integrals (\ref{prod}) we can write
\begin{eqnarray*}
&& \det \Lambda ^{(k)}  =  \frac{1}{2} \left(\frac{n!}{(n-k)!} \frac{m!}{(m-k)! } \right) ^{2}\sum_{ i_{1},\dots, i_{k}, l_{1}, \dots, l_{k}\geq 1}\\
&&\left( \sum_{r=0} ^{(n-k)\wedge (m-k)} r!  C_{ n-k} ^{r} C_{m-k} ^{r} I_{m+n-2k} \left[  f _{i_{1},\dots, i_{k}} \otimes _{r}   g _{l_{1},\dots , l_{k}}  
   -   f _{l_{1},\dots, l_{k}} \otimes _{r}   g _{i_{1},\dots, i_{k}} \right] \right) ^{2}.
\end{eqnarray*}
Taking the mathematical expectation, the isometry of multiple integrals implies that
\begin{eqnarray*}
&& E\det \Lambda ^{(k)}  = \frac{1}{2} \left(\frac{n!}{(n-k)!} \frac{m!}{(m-k)! } \right) ^{2}\sum_{ i_{1},\dots, i_{k}, l_{1}, \dots, l_{k}\geq 1}\\
&& \sum_{r=0} ^{(n-k)\wedge (m-k)} \left(r!  C_{ n-k} ^{r} C_{m-k} ^{r}\right) ^{2} (m+n-2k-2r)!   \\
&& \quad \times 
  \Vert    f _{i_{1},\dots, i_{k}} \tilde{\otimes} _{r}   g _{l_{1},\dots, l_{k}} -  f _{l_{1},\dots, l_{k}} \tilde{\otimes } _{r}  g _{i_{1},\dots, i_{k}}\Vert ^{2}\\
&:=& \sum _{r=0} ^{(n-k)\wedge (m-k)} T_{r}^{(k)},
\end{eqnarray*}
where  
\begin{equation}
\label{tr}
T^{(k)}_{r} =   \frac 12 \alpha_{k,r}  \sum_{ i_{1},\dots, i_{k}, l_{1}, \dots, l_{k}\geq 1}    \Vert    f _{i_{1},\dots, i_{k}} \tilde{\otimes } _{r}   g   _{l_{1},\dots, l_{k}}    
  -   f  _{l_{1}, \dots, l_{k}} \tilde{\otimes} _{r}   g   _{i_{1},\dots, i_{k}} \Vert ^{2} ,
\end{equation}
with
\[
\alpha_{k,r}=  \left( \frac{n!m!}  { (n-k-r)! (m-k-r)!r!} \right)^2 (m+n-2k-2r)!.
\]
We will explicitly  compute the terms $T_{r}^{(k)}$ in (\ref{tr}). To do this, we will need several auxiliary lemmas.
The first one is an immediate consequence of the definition of contraction.

\begin{lemma}\label{l1}
Let $f\in H ^{\odot n} , g\in H^{\odot m}$.  Then for every $k,r\geq 0$ such that $k+r \leq m\wedge n$, 
\begin{equation*}
\sum_{i_{1},\dots, i_{k}\geq 1}   f _{i_{1},\dots,  i_{k}} \otimes _{r}  g _{i_{1}, \dots , i_{k}} = f\otimes _{r+k} g.
\end{equation*}
\end{lemma}
The next lemma summarizes the results in Lemmas 3 and 4 in \cite{T} (see also Lemma 2.2 in \cite{NoRo}). 
\begin{lemma}\label{ll3}
  Assume $f , h \in H^{\odot n}$ and $g, \ell  \in H^{\odot m}$. 
  \begin{enumerate}
  \item[(i)] For every $r=0,\dots , (m-1)\wedge (n-1)$ we have
\begin{equation*}
\langle f   \otimes _{n-r} h,  g \otimes _{m-r}\ell \rangle= \langle f \otimes_{r} g, h \otimes _{r} \ell \rangle.
\end{equation*}
\item[(ii)] The following equality holds
\begin{equation*}
\langle f \tilde{\otimes }g, \ell \tilde{\otimes} h \rangle = \frac{m! n!}{(m+n)! }\sum_{r=0} ^{m\wedge n} C_{n}^{r} C_{m} ^{r} \langle f  \otimes _{r} \ell , h \otimes _{r}g \rangle .
\end{equation*}
\end{enumerate}
\end{lemma}
\vskip0.3cm
We are now ready to calculate the term $T_{0}^{(k)}$. 

\begin{prop}
Let $f\in H ^{\odot n}, g \in H ^{\odot m}$. Let $T_{0}^{(k)}$ be given by (\ref{tr}). Then for every $1\leq k\leq \min (m,n)$
\begin{eqnarray*}
T_{0}^{(k)}&=& \frac{ m! ^{2} n!^{2} } {(m -k)! (n-k) !}  \sum_{s=0}^{ (m-k)\wedge (n-k) } C_{m-k} ^{s} C_{n-k} ^{s}\left[ \Vert f\otimes _{s} g\Vert ^{2}- \Vert f\otimes _{s+k} g\Vert ^{2}\right].
\end{eqnarray*}
\end{prop}
{\bf Proof: }  From (\ref{tr}) we can write
\begin{eqnarray}
  T_{0}^{(k)} &=& \frac{1}{2}  \alpha_{k,0} \sum_{ i_{1},\dots, i_{k}, l_{1}, \dots, l_{k}\geq 1}   
\Vert   f _{i_{1},\dots, i_{k}} \tilde{\otimes}   g _{l_{1},\dots, l_{k}}  -   f _{l_{1},\dots, l_{k}} \tilde{\otimes}    g _{i_{1},\dots, i_{k}} \Vert ^{2}  \nonumber \\  \notag
&   = &\alpha_{k,0} \sum_{ i_{1},\dots, i_{k}, l_{1}, \dots,  l_{k}\geq 1}  
 \Big[ \Vert     f _{i_{1},\dots, i_{k}}   \tilde{\otimes}  g _{l_{1},\dots, l_{k}}  \Vert ^{2}  \\
 &&   - \langle   f  _{i_{1},\dots, i_{k}}   \tilde{\otimes}  g _{l_{1},\dots l_{k}}  ,   g _{i_{1},\dots, i_{k}}   \tilde{\otimes}  f _{l_{1},\dots, l_{k}}  \rangle \Big].    \label{t0-1}
\end{eqnarray}
By Lemma \ref{ll3}, point (ii)  and point (i)
\begin{eqnarray}   \notag
 \Vert     f _{i_{1},\dots, i_{k}}   \tilde{\otimes}  g\ _{l_{1},\dots, l_{k}}  \Vert ^{2} 
& =&  \langle  f _{i_{1},\dots, i_{k}}   \tilde{\otimes}  g _{l_{1},\dots, l_{k}}   , f _{i_{1},\dots, i_{k}}   \tilde{\otimes}  g   _{l_{1},\dots, l_{k}} ) \rangle\\ \notag
 & =&\frac{(m-k)! (n-k)!}{(m+n-2k)!} \sum_{s=0} ^{(m-k)\wedge( n-k)}   C_{m-k}^{s} C_{ n-k}^{s}   \\ \notag
 && \times   \langle  f _{i_{1},\dots, i_{k}}   \otimes _{s} g _{l_{1},\dots, l_{k}}  ,  f _{i_{1},\dots, i_{k}}   \otimes _{s} g _{l_{1},\dots, l_{k}}  \rangle\\   \label{e2}
 &=&\frac{(m-k)! (n-k)!}{(m+n-2k)!} \sum_{s=0} ^{(m-k)\wedge (n-k)} C_{m-k}^{s} C_{ n-k}^{s}  \notag \\
&& \!\!\!\!\!\!\!\!\!\!\!\!\!\! \times   
 \langle  f _{i_{1},\dots, i_{k}}   \otimes _{n-k-s}   f _{i_{1},\dots, i_{k}} ,  g _{l_{1},\dots, l_{k}}    \otimes _{m-k-s} g _{l_{1},\dots, l_{k}}  \rangle.   \label{f1}
\end{eqnarray}
Also, Lemma  \ref{l1}  and Lemma \ref{ll3} point (i)  imply
\begin{eqnarray}
 && \sum_{i_{1},\dots, i_{k}, l_{1}, \dots, l_{k}\geq 1}  \langle  f _{i_{1},\dots, i_{k}}   \otimes _{n-s}   f _{i_{1},\dots, i_{k}} ,  g _{l_{1},\dots, l_{k}}    \otimes _{m-s} g _{l_{1},\dots, l_{k}}  \rangle   \notag  \\
 && \qquad  =\langle  f     \otimes _{n-s}   f  ,  g      \otimes _{m-s} g    \rangle 
 =    \| f\otimes _s g\|^2.  \label{f2}
\end{eqnarray}
On the other hand, using again Lemma \ref{ll3}, point (ii)
\begin{eqnarray}  \notag
&&\langle   f _{i_{1},\dots, i_{k}}   \tilde{\otimes} g _{l_{1},\dots, l_{k}} , g _{i_{1},\dots, i_{k}} ) \tilde{\otimes} f _{l_{1},\dots, l_{k}}  \rangle 
=\frac{(m-k)! (n-k)!}{(m+n-2k)!} \\ 
&& \quad \times  \sum_{s=0} ^{(m-k)\wedge (n-k)} C_{m-k}^{s} C_{ n-k}^{s}    
 \langle  f _{i_{1},\dots, i_{k}}   \otimes _{s} g _{i_{1},\dots, i_{k}}  , f  _{l_{1},\dots, l_{k}}   \otimes _{s} g _{l_{1},\dots, l_{k}}  \rangle.  \label{f3}
\end{eqnarray}
Again, Lemma  \ref{l1}  and Lemma \ref{ll3} point (i)  imply
\begin{eqnarray}
 && \sum_{i_{1},\dots, i_{k}, l_{1}, \dots, l_{k}\geq 1}  \langle  f _{i_{1},\dots, i_{k}}   \otimes _{s} g _{i_{1},\dots, i_{k}}  , f  _{l_{1},\dots, l_{k}}   \otimes _{s} g _{l_{1},\dots, l_{k}}  \rangle   \notag  \\
 && \qquad  =\langle  f     \otimes _{s+k}   g ,  f      \otimes _{s+k} g    \rangle 
 =    \| f\otimes _{s+k} g\|^2.  \label{f4}
\end{eqnarray}
Then, substituting (\ref{f1}), (\ref{f2}), (\ref{f3}) and (\ref{f4}) into (\ref{t0-1}) yields the desired result. 
\qed

It is also possible to compute the terms $T_{r} ^{(k)}$ for every $ 1\leq r \leq (n-k) \wedge (m-k)$ but the  corresponding expressions are more complicated, involving some kind of contractions of contractions. In order to obtain this type of formula we need the following generalization  of  point (ii)  in Lemma \ref{ll3}. 

 For $f , h \in H^{\odot n}$ and $g, \ell  \in H^{\odot m}$ and for $r,s \ge 0$ such that $r+s \le m\wedge n$ 
we denote by  $(f \otimes_r g) \widehat {\otimes}_s(\ell  \otimes_ rh )$  the contraction of $r$ coordinates between $f $ and $g$ and between $\ell $ and $h$, $s$ coordinates between $f $ and $\ell$ and between $g$ and $h$, $n-r-s$ coordinates between $f $ and $h$ and $m-r-s$ coordinates between $g$ an $\ell$. That is,
\begin{eqnarray*}
(f \otimes_r g) \widehat {\otimes}_s(\ell \otimes_ r h)&=&
\sum   \langle  f _{i_1, \dots, i_r, j_1, \dots, j_s, k_1, \dots,  k_{n-r-s} } \rangle     \langle  g_{i_1, \dots, i_r, l_1, \dots, l_s, p_1, \dots,  p_{m-r-s} } \rangle \\
&&\times \langle  \ell_{p_1, \dots, p_r, j_1, \dots, j_s, p_1, \dots,  p_{m-r-s} } \rangle  \langle h_{p_1, \dots, p_r, l_1, \dots, l_s, k_1, \dots,  k_{n-r-s} } \rangle,
\end{eqnarray*}
where the sum runs over all indices greater or equal than one. Notice that
\[
(f \otimes_r g) \widehat {\otimes}_s(\ell  \otimes_ rh)=
(f \otimes_s  \ell ) \widehat {\otimes}_r(g \otimes_ sh).
\]

\begin{lemma}\label{29i-1}
 Assume $f ,h \in H^{\odot n}$ and $g,\ell \in H^{\odot m}$.    Then for every $r=0,\dots, (m-1)\wedge (n-1)$ we have
\begin{equation*}
\langle f  \tilde{\otimes }_{r}g,  \ell  \tilde{\otimes}_{r} h \rangle =\frac{(n-r)! (m-r)!}{(m+n-2r)!}  \sum_{s=0} ^{(m-r) \wedge (n-r)} C_{n-r}^{s} C_{m-r} ^{s} (f \otimes_r g) \widehat {\otimes}_s(\ell  \otimes_ rh).
\end{equation*}
\end{lemma}
{\bf Proof: }  We can write
\begin{equation*}
 f  \tilde{\otimes }_{r}g= \sum_{i_1, \dots, i_r}    f _{i_1, \dots, i_r}   \tilde{\otimes}   g_{i_1, \dots, i_r}   
\end{equation*}
and
\begin{equation*}
\ell  \tilde{\otimes }_{r}h= \sum_{i_1, \dots, i_r}  \ell_{i_1, \dots, i_r}   \tilde{\otimes}  h_{i_1, \dots, i_r} .
\end{equation*}
Then,  Lemma \ref{ll3} point  (ii) gives
\begin{eqnarray*}
&&\langle f  \tilde{\otimes }_{r}g, \ell \tilde{\otimes}_{r} h \rangle 
 =  \sum_{i_1, \dots, i_r, l_1, \dots l_r}   \langle   f _{i_1, \dots, i_r}   \tilde{\otimes}  g_{i_1, \dots, i_r}   ,
\ell _{l_1, \dots, l_r}   \tilde{\otimes}  h_{l_1, \dots, l_r}\rangle   \\
&=& \frac{(n-r)! (m-r)!}{(m+n-2r)!} 
\sum_{s=0} ^{(m-r) \wedge (n-r)} C_{n-r}^{s} C_{m-r} ^{s}   \\
&& \langle   f  _{i_1, \dots, i_r} \otimes _s \ell _{l_1, \dots, l_r}  , h_{l_1, \dots, l_r} \otimes_s
g_{i_1, \dots, i_r}
 \rangle,
\end{eqnarray*}
which implies the desired result.
\qed

Notice that for $r=0$,
\[
(f \otimes g) \widehat {\otimes}_s(\ell  \otimes h) = \langle f \otimes_s \ell , h\otimes_s g  \rangle,
\]
so Lemma \ref{ll3} point (ii)  is a particular case of Lemma  \ref{29i-1} when $r=0$.

\begin{prop}
Let $(F,G)= (I_{n}(f), I_{m}(g))$ with $f\in H^{\odot n}$ and $g\in H^{\odot m}$.   Then, for  every $r=1,\dots, (n-k)\wedge(m-k)$
\begin{eqnarray}  \notag
T_{r} ^{(k)}& =& \beta_{k,r}\sum_{s=0} ^{(n-k-r)\wedge (m-k-r)}  C_{n-k-r}^{s}C_{m-k-r}^{s} \\
&& \times \left( 
(f\otimes_r g)\widehat{\otimes}_s  (g\otimes _rf)- (f\otimes_r g)\widehat{\otimes}_{s+k}  (g\otimes _rf) \right), \label{tr1}
\end{eqnarray}
where
\[
\beta_{k,r}= \frac {n! ^{2}m! ^{2}}  {(n-k-r)! (m-k-r)! (r!)^2}.
\]
\end{prop}

\noindent
{\bf Proof: }   From (\ref{tr}) we can write
\begin{eqnarray}
T_{r}^{(k)}&=&  \alpha_{k,r} \sum_{i_{1},\dots,i_{k}, l_{1},\dots, l_{k}\geq 1}
\left[ \Vert  f  _{i_{1},\dots , i_{k}}  \tilde{\otimes} _{r}   g _{l_{1},\dots,l_{k}}) \Vert ^{2}   \right. \nonumber \\
&& \left. 
 - \langle  f  _{i_{1},\dots, i_{k}}  \tilde{\otimes}_{r}  g _{l_{1},\dots,l_{k}}  ,  f _{l_{1},\dots,l_{k}} )\tilde{\otimes} _{r}   g _{i_{1},\dots,i_{k}}  \rangle \right].\label{28a3}
\end{eqnarray} 
Applying Lemma \ref{29i-1} yields
\begin{eqnarray}  \notag
&& \Vert  f  _{i_{1},\dots, i_{k}}  \tilde{\otimes} _{r}   g _{l_{1},\dots,l_{k}}  \Vert ^{2} 
  \langle  f  _{i_{1},\dots, i_{k}}  \tilde{\otimes} _{r}   g _{l_{1},\dots,l_{k}}  ,  f  _{i_{1},\dots, i_{k}} )\tilde{\otimes} _{r}   g _{l_{1},\dots,l_{k}}  \rangle\\  \notag
&=&\frac{(n-k-r)! (m-k-r)!}{(m+n-2k-2r)!} \sum_{s=0} ^{(n-k-r)\wedge(m-k-r)}  C_{n-k-r}^{s}C_{m-k-r}^{s}   \\
&&\times
\left(   f _{i_{1},\dots, i_{k}}  \otimes _{r}  g _{l_{1},\dots,l_{k}}   \right) \widehat{\otimes}_s \left(  g _{l_{1},\dots,l_{k}}  \otimes_r f  _{i_{1},\dots, i_{k}}     \right).  \label{x1}
\end{eqnarray}
Notice that 
\begin{equation}  \label{x2}
\sum_{i_{1},\dots,i_{k}, l_{1},\dots, l_{k}\geq 1}\left(   f  _{i_{1},\dots, i_{k}}  \otimes _{r}  g\ _{l_{1},\dots,l_{k}}   \right)  
  \widehat{\otimes}_s \left(   g _{l_{1},\dots,l_{k}}  \otimes_r f  _{i_{1},\dots, i_{k}}     \right)
= (f\otimes_r g)\widehat{\otimes}_s  (g\otimes _rf).
\end{equation} 
Analogously, we get
\begin{eqnarray}  \notag
&& \langle  f  _{i_{1},\dots, i_{k}}  \tilde{\otimes} _{r}   g _{l_{1},\dots,l_{k}}  ,  f  _{l_{1},\dots, l_{k}}  \tilde{\otimes} _{r}   g _{i_{1},\dots,i_{k}}  \rangle\\  \notag
&=&\frac{(n-k-r)! (m-k-r)!}{(m+n-2k-2r)!} \sum_{s=0} ^{(n-k-r)\wedge(m-k-r)}  C_{n-k-r}^{s}C_{m-k-r}^{s}   \\
&&\times
\left(   f _{i_{1},\dots, i_{k}}  \otimes _{r}   g _{l_{1},\dots,l_{k}}   \right) \widehat{\otimes}_s \left(   g _{i_{1},\dots,i_{k}}) \otimes_r f  _{l_{1},\dots, l_{k}}    \right),  \label{x3}
\end{eqnarray}
and
\begin{equation}  \label{x4}
\sum_{i_{1},\dots,i_{k}, l_{1},\dots, l_{k}\geq 1}\left(   f   _{i_{1},\dots, i_{k}}  \otimes _{r}  g _{l_{1},\dots,l_{k}}   \right)  
 \widehat{\otimes}_s \left(   g\ _{i_{1},\dots,i_{k}}  \otimes_r f   _{l_{1},\dots, l_{k}}     \right)
= (f\otimes_r g)\widehat{\otimes}_{s+k}  (g\otimes _rf).
\end{equation} 
Substituting (\ref{x1}), (\ref{x2}), (\ref{x3}) and (\ref{x4}) into (\ref{28a3}) we obtain the desired formula. 
\qed

In the particular case $n=m$,  the expression (\ref{tr1}) can be written as
\begin{eqnarray}\label{tr2}
T_{r} ^{(k)} =  \sum_{s=0} ^{n-k-r}     T^{(k)}_{r,s},
\end{eqnarray}
where
\[
T^{(k)}_{r,s}= \frac { (n!)^4}  {((n-k-r)! r!)^2} (C_{n-k-r}^{s})^2 \left( 
(f\otimes_r g)\widehat{\otimes}_s  (g\otimes _rf)- (f\otimes_r g)\widehat{\otimes}_{s+k}  (g\otimes _rf) \right).
\]
The  last term in (\ref{tr2}) obtained for $r= n-k$ is given by the following expression.

\begin{corollary}  \label{col1}
Let $(F,G)= (I_{n}(f), I_{n}(g))$ with $f, g\in H^{\odot n}$. Then for $k=1,\dots,n-1$
\begin{equation*}
T_{n-k} ^{(k)}=\frac{n! ^{4}}{(n-k)! ^{2}} \left[ \Vert f\otimes _{n-k} g\Vert ^{2}- \langle f\otimes _{n-k} g, g\otimes _{n-k}f\rangle \right].
\end{equation*}

\end{corollary} 
{\bf Proof: } When $r=n-k$, there is only one terms in the sum (\ref{tr2}), obtained for $s=0$. It is easy to see that, 
\[
(f\otimes_{n-k} g)\widehat{\otimes}_0  (g\otimes _{n-k}f)= (f\otimes g)\widehat{\otimes}_{n-k}  (g\otimes  f)=\Vert f\otimes _{n-k} g\Vert ^{2}
\]
and
 \[
(f\otimes_{n-k} g)\widehat{\otimes}_{k}  (g\otimes _{n-k})  =\langle f\otimes _{n-k} g, g\otimes _{n-k}f\rangle.
 \]
 \qed

We obtain the following expression for the determinant of the $k$th Malliavin matrix.

\begin{theorem}\label{t1}
Let $f\in H ^{\odot n}, g \in H ^{\odot m} $. Then for every $1\le k\le m\wedge n$,
\begin{eqnarray*}
E \det \Lambda ^{(k)} &=&  \frac{ m! ^{2} n!^{2} } {(m-k)! (n-k) !}  \sum_{s=0}^{ (m-k) \wedge (n-k) } C_{m-k} ^{s} C_{n-k} ^{s}\\
&&\times \left[ \Vert f\otimes _{s} g\Vert ^{2}- \Vert f\otimes _{s+k} g\Vert ^{2}\right] + R_{m,n,k}  ,
\end{eqnarray*}
where  $R_{m,n,k}=\sum_{r=1} ^{(m-k)\wedge (n-k)}
T^{(k)}_{r}$ and  $T^{(k)}_{r}$ is given by
 (\ref{tr1}). 
\end{theorem}

In the case of multiple integrals of the same order (i.e. $m=n$) we have the following result.

\begin{corollary}  \label{col2}
If $f,g \in H ^{\odot n} $, the determinant of the $k$th iterated Malliavin matrix of $(F,G)= (I_{n}(f), I_{n}(g))$ can be written as 
\[
E \det \Lambda ^{(k)} =\frac{ n! ^{4}}{(n-k) ! ^{2}} \sum_{s=0} ^{n-k}  (C_{n-k}^{s}) ^{2}\left( \Vert f\otimes _{s}g\Vert ^{2} -\Vert f\otimes _{s+k} g\Vert ^{2}\right) + R_{n,n,k}.
\]
\end{corollary}

\begin{example}
Suppose $m=n=3$ and $k=2$. Then
\[ 
E\det \Lambda ^{(2)}= (3!) ^{4} \left[  \Vert f\otimes _{0} g \Vert ^{2} - \Vert f\otimes _{2} g \Vert ^{2} 
+ \Vert f\otimes _{1} g\Vert ^{2}- \Vert f\otimes _{3} g \Vert ^ {2}  \right]+R_{3,3,2}.
\]
Suppose $m=n=4$ and $k=2$. Then
\[
E\det \Lambda ^{(2)}= \frac{(4!) ^{4} }{ 2! ^{2}} \left[ \Vert f\otimes _{0} g \Vert ^{2} - \Vert f\otimes _{2} g \Vert ^{2}  + 4( \Vert f\otimes _{1}g\Vert ^{2} - \Vert f\otimes _{3} g\Vert  ^{2} )\right]+ R_{4,4,2}.
\]
\end{example}

Our next objective is to relate the expectation of the iterated Malliavin matrix  $ E \det \Lambda ^{(s)}$ with the covariance matrix of the vector $(F,G)$ in the case $n=m$. We recall that
\[
\det C= n!^2 [\|f\|^2 \|g\|^2- \langle f,g \rangle^2].
\]
 \begin{theorem} \label{thm1}
 For any  $f,g \in H ^{\odot n} $, if $F=I_n(f)$ and $G=I_n(g)$, we have
 \[
  \sum_{s=2} ^{\left[ \frac{n-1}2 \right]}    \frac {n(n-2s)}{s!^2} E\det \Lambda^{(s)} +
    (n-1)^2 E \det \Lambda ^{(1)} \ge n^2 \det C.
 \]
 \end{theorem}
 {\bf Proof: } 
 From Corollary \ref{col2}, taking into account that $R_{n,n,1} \ge 0$, we can write
 \begin{eqnarray*}
 E \det \Lambda ^{(1)} &\ge&  [nn!]^2
\sum_{s=0} ^{n-1}  (C_{n-1}^{s}) ^{2}\left( \Vert f\otimes _{s}g\Vert ^{2} -\Vert f\otimes _{s+1} g\Vert ^{2}\right) \\
&= & n^2 \det C +[nn!]^2 \sum_{s=1} ^{n-1}   \left( (C_{n-1}^{s}) ^{2}- (C_{n-1}^{s-1})^2 \right)  \Vert f\otimes _{s}g\Vert ^{2}.
\end{eqnarray*}
Notice that $ (C_{n-1}^{s}) ^{2}- (C_{n-1}^{s-1})^2 =-[(C_{n-1}^{n-s}) ^{2}- (C_{n-1}^{n-1-s})^2]$. Therefore, we conclude that
 \begin{eqnarray}  \notag
 E \det \Lambda ^{(1)} &\ge&   n^2 \det C +[nn!]^2 \sum_{s=1} ^{\left[ \frac{n-1}2 \right]}   \left( (C_{n-1}^{s}) ^{2}- (C_{n-1}^{s-1})^2 \right)  \\  \notag
 &&\times  \left( \Vert f\otimes _{s}g\Vert ^{2}- \Vert f\otimes _{n-s}g\Vert ^{2} \right)\\
 &=&  n^2 \det C + \sum_{s=1} ^{\left[ \frac{n-1}2 \right]}    \gamma_{n,s}  \left( \Vert f\otimes _{s}g\Vert ^{2}- \Vert f\otimes _{n-s}g\Vert ^{2} \right),  \label{h2}
\end{eqnarray}
where
\[
\gamma_{n,s}= \left(\frac { n!^2}{(n-s)! s!} \right)^2 n(n-2s).
\]
Notice that $\gamma_{n,s} \ge 0$ if $s \le \left[ \frac{n-1}2 \right]$. We can write, using Lemma  \ref{ll3} point (i) and Corollary \ref{col1}
\begin{eqnarray}  \notag
\Vert f\otimes _{s}g\Vert ^{2}- \Vert f\otimes _{n-s}g\Vert ^{2} &=&  \Vert f\otimes _{s}g\Vert ^{2}- 
\langle f\otimes _s g, g\otimes _s f \rangle  \\  \notag
&& - \left( \Vert f\otimes _{n-s}g\Vert ^{2} -\langle f\otimes _{n-s} g, g\otimes _{n-s} f \rangle \right)\\
&\ge & - \frac {(n-s)!^2} { n!^4}T^{(s)}_{n-s}   \ge- \frac {(n-s)!^2} { n!^4}  E\det \Lambda^{(s)}.  \label{h1}
\end{eqnarray}
Substituting (\ref{h1}) into  (\ref{h2}) yields
\[
 E \det \Lambda ^{(1)} \ge n^2 \det C - \sum_{s=1} ^{\left[ \frac{n-1}2 \right]}   \frac {n (n-2s)}{s!^2} E\det \Lambda^{(s)},
 \]
which implies the desired result.
\qed

\begin{remark}
In the particular case $n=2$ we obtain $E \det \Lambda ^{(1)} \ge 4 \det C$, which was proved in \cite{NoNuPo}.
 For $n=3$ we get  $E \det \Lambda ^{(1)} \ge \frac 94 \det C$, and for $n=4$,
 $E \det \Lambda ^{(1)} \ge \frac {16}9 \det C$. Only if $n\ge 5$ we need the expectation of the iterated Malliavin matrix to control the determinant of the covariance matrix.
\end{remark}

\section{The density of a couple of multiple integrals}

In this section, we    show that a random vector of dimension  2 whose components are multiple integrals in the same Wiener chaos either admits a density with respect to the Lebesque measure, or its components are proportional. We also show that a necessary and sufficient condition for such a vector to not have a density is that at least one of its iterated Malliavin matrices vanishes almost surely.
 In the sequel we fix a vector $(F,G)= (I_{n}(f), I_{n}(g))$ with $f, g \in H^{\odot n}$.

\vskip0.2cm
In the following result we show that, if the determinant of an iterated  Malliavin matrix of a couple of multiple integrals  vanishes, the determinant of the any other iterated  Malliavin matrices will vanish.

\begin{prop}\label{1-k} Let $1\leq k,l\leq n$ with $k\not= l$. Then $E\det \Lambda ^{(k)} =0$  if and only if 
$E\det \Lambda ^{(l)} =0$.

\end{prop}
{\bf Proof: }  Assume first that $k=1$ and $l=2$. Suppose that $E\det \Lambda ^{(1)} =0$ and let us prove that $E\det \Lambda ^{(2)} =0$.
 Since $\det \Lambda ^{(1)}=0$ a.s., from (\ref{e1}) we obtain
\begin{equation} \label{w10}
D_{j}F D_{i}G= D_{i}F D_{j} G   \quad  \mbox {a.s.}
\end{equation}
for any $i,j\ge 1$ (recall that $D_jF =  DF \otimes_1 e_j$). That is,
\begin{equation} \label{w11} 
DF D_iG = DG D_iF   \quad  \mbox {a.s.},
 \end{equation}
for any $i\ge 1$.
Let us apply the divergence operator $\delta$ (the adjoint of $D$) to both members of equation (\ref{w11}). 
From  (\ref{delta}) we obtain $\delta DF= nF$ and $\delta DG =nG$. Using Proposition 1.3.3 in \cite{N}, we get 
\[
nFD_iG - \langle DF, DD_iG\rangle_H=nGD_iF - \langle DG, DD_iF\rangle_H  \quad  \mbox {a.s.},
\]
which can be written as (using the notation (\ref{fjj}))
\begin{eqnarray*}
&&I_{n}(f) I_{n-1} (g_i)- (n-1)  \sum_{j=1}^\infty I_{n-2} (g_{ij}) I_{n-1}(f_j)\\
&&=I_{n}(g) I_{n-1} (f_i)- (n-1) \sum_{j=1}^\infty I_{n-2} (f_{ij}) I_{n-1}(g_j)  \quad  \mbox {a.s.}
\end{eqnarray*}
 By the product formula (\ref{prod}), the above relation becomes
\begin{eqnarray*}
&&I_{2n-1} (f\tilde{\otimes} g_i)+ \sum _{k=1} ^{n-1} (k! C_{n}^{k}C_{n-1} ^{k}-(n-1) (k-1)! C_{n-1} ^{k-1} C_{n-2} ^{k-1} ) \\
&&\times I_{2n-1-2k} (f\tilde{\otimes} _{k} g_i) \\
&&=I_{2n-1} (g\tilde{\otimes }f_i)+ \sum _{k=1} ^{n-1} (k! C_{n}^{k}C_{n-1} ^{k}- (n-1)(k-1)! C_{n-1} ^{k-1} C_{n-2} ^{k-1} )  \\
&&\times I_{2n-1-2k} (g\tilde{\otimes} _{k} f_i)  \quad  \mbox {a.s.}
\end{eqnarray*}
 By identifying the terms in each Wiener chaos, we obtain 
 \[
 f\tilde{\otimes}_k g_i = g\tilde{\otimes}_k f_i
 \]
 for any $i\ge 1$ and for any $k=0,\dots, n-1$. A further application of the product formula for multiple integrals yields 
\begin{equation*}
FDG = GDF  \quad  \mbox {a.s.}
\end{equation*}
We differentiate the above relation in the Malliavin sense and we have
\begin{equation*}
FD^{(2)}_{ij} G + D_{i}FD_{j} G = GD^{(2)} _{ij}F + D_{i}G D_{j} F   \quad  \mbox {a.s.}
\end{equation*}
for every $i,j\ge 1$. By (\ref{w10}),
\begin{equation*}
F D^{(2)}G = GD ^{(2)} F   \quad  \mbox {a.s.}
\end{equation*}
and this clearly implies that $\det \Lambda ^{(2) }=0$ a.s.

Suppose now that $E\det \Lambda ^{(2)}=0$.  Then  $\Lambda ^{(2)}=0$ a.s. and from (\ref{e1}) we get
\[
D^{(2)}_{ij} F D^{(2)}_{pq} G = D^{(2)}_{pq}F D^{(2)}_{ij}G   \quad  \mbox {a.s.}
\]
for any $i,j,p,q \ge 1$. This implies
\begin{equation} \label{w12}
D D_{i} F D^{(2)}_{pq} G =  D D_{i}G D^{(2)}_{pq}F    \quad  \mbox {a.s.}
\end{equation}
for any $i,p,q \ge 1$ Applying 
  the divergence operator $\delta $  to equation (\ref{w12}) yields
  \[
  (n-1)D_{i} F D^{(2)}_{pq} G- \langle  D D_{i} F ,DD^{(2)}_{pq} G\rangle_H=  (n-1) D_{i}G D^{(2)}_{pq}F  -
  \langle D D_{i}G , DD^{(2)}_{pq}F  \rangle_H 
  \]
  a.s. This equality can be written as
  \begin{eqnarray*}
&&I_{n-1}(f_i) I_{n-2} (g_{pq})- (n-2)  \sum_{j=1}^\infty I_{n-3} (g_{pq}) I_{n-2}(f_{ij})\\
&&=I_{n}(g_i) I_{n-1} (f_{pq})- (n-2) \sum_{j=1}^\infty I_{n-3} (f_{hl}) I_{n-2}(g_{pq})  \quad  \mbox {a.s.}
\end{eqnarray*}
  By  the product formula  for multiple integrals  we get for every $j,p,q\ge 1$
\begin{eqnarray*}
&&I_{2n-3}\left( f_i\tilde{\otimes} g_{pq}\right)  +\sum_{k=1}^{n-2} \left[ k! C_{n-2}^{k} C_{n-1}^{k}- (n-2)(k-1)! C_{n-2}^{k-1} C_{n-3}^{k-1}\right]  \\
&&\qquad \times I_{2n-3-2k} \left(f_i\tilde{\otimes}_{k} g_{pq} \right)\\
&=& I_{2n-3}\left( g_i\tilde{\otimes f } _{pq}\right)  +\sum_{k=1}^{n-2} \left[ k! C_{n-2}^{k} C_{n-1}^{k}- (n-2)(k-1)! C_{n-2}^{k-1} C_{n-3}^{k-1}\right] \\
&&\qquad  \times  I_{2n-3-2k} \left( g_i \tilde{\otimes}_{k} f_{pq}\right) \quad  \mbox{ a.s. } 
\end{eqnarray*}
Identifying the coefficients of each Wiener chaos   we obtain 
\[
f_i\tilde{\otimes}_k g_{pq} = g_i \tilde{\otimes}_k f_{pq}
\]
for any $i,p,q \ge 1$ and for any $k=0, \dots, n-2$.
This implies
\begin{equation} \label{w13}
f_i\tilde{\otimes}_k g_{q} = g_i\tilde{ \otimes}_k f_{q}
\end{equation}
for any $i,q \ge 1$ and for any $k=0, \dots, n-1$. Applying again the product formula for multiple integrals
(\ref{w13}) leads to
\[
D_i F D_q G = D_i G D_qF  \quad  \mbox{ a.s. } ,
\]
for any $i,q\ge 1$, which implies 
$\det \Lambda ^{(1)} =0$ a.s.
By iterating the above argument, we easily find that $\det \Lambda ^{(k)}= 0$  a.s. is equivalent to  $\det \Lambda ^{(l)}= 0$  a.s., for every $1\leq k,l\leq n$ with $k\not=l$. \qed

\begin{corollary}
The vector $(F,G)= (I_{n}(f), I_{n}(g))$ does not admit a density if and only if there exists $k\in \{1,\dots, n\}$ such that $E\det \Lambda ^{(k)}= 0$.
\end{corollary}
{\bf Proof: } It is a consequence of Proposition \ref{1-k} and of Theorem 3.1 in \cite{NoNuPo}. \qed

\begin{theorem}

Let  $f,g \in H^{\odot n}$ be symmetric tensors. Then the random vector $(F,G)= (I_{n}(f), I_{n} (g))$ does not admit a density if and only if 
$\det C =0$
where $C$ denotes the covariance matrix of $(F,G)$.
 In other words, the vector $(F, G)$ does not admit a density if and only if its components are proportional.

\end{theorem}
{\bf Proof: }  If  $\det C=0$, the random variables $F$ and $G$ are proportional and the law of $(F,G)$ is not absolutely continuous with respect to the Lebesgue measure.
Suppose that the law of the random vector $(F,G)$ is not absolutely continuous with respect to the Lebesque measure. Then, from the results of \cite{NoNuPo} we know that  $E\det \Lambda ^{(1)}$=0. By Proposition \ref{1-k}, $E\det \Lambda ^{(k)}=0$ for $k=1,\dots,n$.
Then Theorem \ref{thm1} implies
$\det C=0$ (notice also that $\det C=0$ because  $C=n! \Lambda^{(n)}$).\qed

\end{document}